\documentclass{article}

\PassOptionsToPackage{numbers, compress}{natbib}

\usepackage[preprint]{nips_2018}

\usepackage[utf8]{inputenc} 
\usepackage[T1]{fontenc}    
\usepackage{hyperref}       
\usepackage{url}            
\usepackage{booktabs}       
\usepackage{amsfonts}       
\usepackage{nicefrac}       
\usepackage{microtype}      

\usepackage{times}
\usepackage{epsfig}
\usepackage{graphicx}
\usepackage{amsmath}
\usepackage{amssymb}

\usepackage{arydshln}
\usepackage{pdftricks}
\usepackage{epsfig}
\usepackage{pst-grad} 
\usepackage{pst-plot} 
\usepackage[space]{grffile} 
\usepackage{etoolbox} 
\makeatletter 
\patchcmd\Gread@eps{\@inputcheck#1 }{\@inputcheck"#1"\relax}{}{}
\makeatother

\title{Rethinking floating point for deep learning}

\author{
Jeff Johnson\\
Facebook AI Research\\
New York, NY\\
{\tt\small jhj@fb.com}
}

\begin{document}

\maketitle

\begin{abstract}
  Reducing hardware overhead of neural networks for faster or lower power inference and training is an active area of research. 
  Uniform quantization using integer multiply-add has been thoroughly investigated, which requires learning many quantization parameters, fine-tuning training or other prerequisites.
  Little effort is made to improve floating point relative to this baseline; it remains energy inefficient, and word size reduction yields drastic loss in needed dynamic range.
  We improve floating point to be more energy efficient than equivalent bit width integer hardware on a 28 nm ASIC process while retaining accuracy in 8 bits with a novel hybrid log multiply/linear add, Kulisch accumulation and tapered encodings from Gustafson's posit format.
With no network retraining, and drop-in replacement of all math and float32 parameters via round-to-nearest-even only, this open-sourced 8-bit log float is within 0.9\% top-1 and 0.2\% top-5 accuracy of the original float32 ResNet-50 CNN model on ImageNet. Unlike int8 quantization, it is still a general purpose floating point arithmetic, interpretable out-of-the-box.
Our 8/38-bit log float multiply-add is synthesized and power profiled at 28 nm at $0.96\times$ the power and $1.12\times$ the area of 8/32-bit integer multiply-add. In 16 bits, our log float multiply-add is 0.59$\times$ the power and 0.68$\times$ the area of IEEE 754 float16 fused multiply-add, maintaining the same signficand precision and dynamic range, proving useful for training ASICs as well.
\end{abstract}

\section{Introduction}\label{sec:introduction}

Reducing the computational complexity of neural networks (NNs) while maintaining accuracy encompasses a long line of research in NN design, training and inference.
Different computer arithmetic primitives have been considered, including fixed-point~\cite{Lin2016}, uniform quantization via 8 bit integer~\cite{Jacob2018}, ternary~\cite{DBLP:journals/corr/LiL16} and binary/low-bit representations~\cite{rastegari2016xnor, courbariaux2016binarized, cai2017deep}. Some implementations are efficiently implemented on CPU/GPU ISAs~\cite{vanhoucke2011improving, DBLP:journals/corr/abs-1712-02427}, while others demand custom hardware~\cite{gupta2015deep}.
Instead of developing quantization techniques increasingly divorced from the original implementation, we seek to improve floating point itself, and let word size reduction yield efficiency for us. It is historically known to be up to 10$\times$ less energy efficient in hardware implementations than integer math~\cite{horowitz20141}. Typical implementation is encumbered with IEEE 754 standard compliance~\cite{zuras2008ieee}, demanding specific forms such as fused multiply-add (FMA) that we will show as being inefficient and imprecise. Memory movement (SRAM/DRAM/flip-flops) dominates power consumption; word bit length reduction thus provides obvious advantages beyond just reducing adder and multiplier area.



We explore encodings to better capture dynamic range with acceptable precision in smaller word sizes, and more efficient summation and multiplication (Sections~\ref{sec:encoding}-\ref{sec:mul}), for a reduction in chip power and area. Significant inspiration for our work is found in logarithmic number systems (LNS)~\cite{coleman2000arithmetic} and the work of Miyashita et al.~\cite{miyashita2016convolutional} that finds logarithmic quantizers better suited to data distributions in NNs, and alternative visions of floating point from Gustafson~\cite{Gustafson2015, Gustafson2017} and Kulisch~\cite{Kulisch2012}. We sidestep prior LNS design issues with numerical approximation and repurpose ideas from Gustafson and Kulisch, producing a general-purpose arithmetic that is effective on CNNs~\cite{he2016deep} without quantization tinkering or re-training (Section~\ref{sec:analysis}), and can be as efficient as integer math in hardware (Section~\ref{sec:ASIC}).

\section {Floating point variants for NNs}\label{sec:prior}

There are few studies on NNs for floating point variants beyond those provided for in CPU/GPU ISAs. \cite{dettmers20158} shows a kind of 8 bit floating point for communicating gradients, but this is not used for general computation. Flexpoint~\cite{koster2017flexpoint} and the Brainwave NPU~\cite{fowers2018configurable} use variants of \textit{block floating point}~\cite{wilkinson1963rounding}, representing data as a collection of significands with a shared exponent. This requires controlled dynamic range variation and increased management cost, but saves on data movement and hardware resources. For going to 8 bits in our work, we seek to improve the encoding and hardware for a reasonable tradeoff between dynamic range and precision, with less machinery needed in software.

For different precisions, \cite{DiCecco2017} shows reduced-precision floating point for training smaller networks on MNIST and CIFAR-10, with (6, 5)\footnote{Throughout, $(e, s)$-float refers to IEEE 754-style floating point, with sign bit, e-bit biased exponent and s-bit 0.s or 1.s fixed point significand; float16/float32 are shorthand for IEEE 754 binary16/binary32.} floating point without denormal significands being comparable to float32 on these examples. (8, 7) \textit{bfloat16} is available on Google's TPUv2~\cite{bfloat16}. This form maintains the same normalized exponent range as float32, except with reduced precision and smaller multipliers. However, the forms of encoding and computation for many of these variants are not substantially different than implementations available with common ISAs, hardened FPGA IP, and the like. We will seek to improve the encoding, precision and computation efficiency of floating point to find a solution that is quite different in practice than standard $(e, s)$ floating point.

\begin{table}
{
\caption{\label{tab:range}
Dynamic range and significand fractional precision of math types considered
}}
\centering \begin{tabular}{llll}
\textbf{Word} & \textbf{Encoding} & \textbf{Range in decibels} & \textbf{Fraction} \\
\textbf{bits} & \textbf{type} & $20 \log_{10}(f_{max}/f_{min})$ & \textbf{bits (max)}\\
\hline
8 & symmetric integer $[-2^{7}+1, 2^7 - 1]$ & 42.1 & ---\\
8 & (8, 0) posit or (8, 0, $\alpha$, $\beta$, $\gamma$) log& 72.2 & 5 \\
8 & (4, 3) float (w/o denormals) & 83.7 & 3\\
16 & symmetric integer $[-2^{15}+1, 2^{15}-1]$ & 90.3 & ---\\
8 & (4, 3) float (w/ denormals) & 101.8 & 3\\
8 & (8, 1) posit or (8, 1, $\alpha$, $\beta$, $\gamma$) log& 144.5 & 4\\
16 & (5, 10) float16 (w/o denormals) & 180.6 & 10\\
16 & (5, 10) float16 (w/ denormals) & 240.8 & 10\\
12 & (12, 1) posit or (12, 1, $\alpha$, $\beta$, $\gamma$) log& 240.8 & 8\\
8  & (8, 2) posit or (8, 2, $\alpha$, $\beta$, $\gamma$) log& 289.0 & 3\\
16 & (16, 1) posit or (16, 1, $\alpha$, $\beta$, $\gamma$) log& 337.2 & 12\\
\end{tabular}
\end{table}

\section {Space-efficient encodings}\label{sec:encoding}

IEEE 754-style fixed width field encodings are not optimal for most data distributions seen in practice; float32 maintains the same significand precision at $10^{-10}$ as at $10^{10}$. 
Straightforward implementation of this design in 8 bits will result in sizable space encoding NaNs, $\sim 6\%$ for (4, 3) float. Denormals use similar space and are expensive in hardware~\cite{muller2010handbook}; not implementing them restricts the dynamic range of the type (Table~\ref{tab:range}).
\textit{Tapered floating point} can solve this problem: within a fixed-sized word, exponent and significand field size varies, with a third field indicating relative size. To quote Morris (1971): ``users of floating-point numbers are seldom, if ever, concerned \textit{simultaneously} with loss of accuracy and with overflow. If this is so, then the range of possible representation can be extended [with tapering] to an extreme degree and the slight loss of accuracy will be unnoticed.''~\cite{morris1971tapered}



A more efficient representation for tapered floating point is the recent \textit{posit} format by Gustafson \cite{Gustafson2017}. It has no explicit size field; the exponent is encoded using a Golomb-Rice prefix-free code~\cite{golomb1966run, Lindstrom2018}, with the exponent $e$ encoded as a Golomb-Rice quotient and remainder $(q, r)$ with $q$ in unary and $r$ in binary (in posit terminology, $q$ is the \textit{regime}). Remainder encoding size is defined by the \textit{exponent scale} $s$, where $2^s$ is the Golomb-Rice divisor. Any space not used by the exponent encoding is used by the significand, which unlike IEEE 754 always has a leading 1; gradual underflow (and overflow) is handled by tapering. A posit number system is characterized by $(N, s)$, where $N$ is the word length in bits and $s$ is the exponent scale. The minimum and maximum positive finite numbers in $(N, s)$ are $f_{min} = 2^{-(N-2)2^{s}}$ and $f_{max} = 2^{(N-2)2^{s}}$.  The number line is represented much as the projective reals, with a single point at $\pm \infty$ bounding $-f_{max}$ and $f_{max}$. $\pm \infty$ and 0 have special encodings; there is no NaN. The number system allows any choice of $N \geq 3$ and $0 \leq s \leq N-3$. 

$s$ controls the dynamic range achievable; e.g., 8-bit $(8, 5)$-posit $f_{max} = 2^{192}$ is larger than $f_{max}$ in float32. $(8, 0)$ and $(8, 1)$ are more reasonable values to choose for 8-bit floating point representations, with $f_{max}$ of 64 and 4096 accordingly. Precision is maximized in the range $\pm [2^{-(s+1)}, 2^{s+1})$ with $N-3-s$ significand fraction bits, tapering to no fraction bits at $\pm f_{max}$.

\section {Accumulator efficiency and precision}\label{sec:accum}

A sum of scalar products $\sum_i a_i b_i$ is a frequent operation in linear algebra. For CNNs like ResNet-50~\cite{he2016deep}, we accumulate up to 4,608 (2d convolution with $k = 3 \times 3, c_{in} = 512$) such products.


Integer addition is associative (excepting overflow); the order of operations does not matter and thus it allows for error-free parallelization. In typical accelerator use, the accumulation type is 32 bits.
Typical floating point addition is notorious for its lack of associativity; this presents problems with reproducibility, parallelization and rounding error~\cite{muller2010handbook}. Facilities such as \textit{fused multiply-add} (FMA) that perform a sum and product $c + a_ib_i$ with a single rounding can reduce error and further pipeline operations when computing sums of products. Such machinery cannot avoid rounding error involved with tiny (8-bit) floating point types; the accumulator can become larger in magnitude than the product being accumulated into it, and the significand words no longer overlap as needed even with rounding (yielding $c + ab = c$); increasing accumulator size a bit only defers this problem.


There is a more efficient and precise method than FMA available. A \textit{Kulisch accumulator}~\cite{Kulisch2012} is a fixed point register that is wide enough to contain both the largest and smallest possible scalar product of floating point values $\pm (f^2_{max} + f^2_{min})$. It provides associative, error-free calculation (excepting a single, final rounding) of a sum of scalar floating point products; a float significand to be accumulated is shifted based on exponent to align with the accumulator for the sum. Final rounding to floating point is performed after all sums are made. A similar operation known as \textit{Auflaufenlassen} was available in Konrad Zuse's Z3 as early as 1941~\cite{Kulisch2002}, though it is not found in modern computers.

We will term this operation of summing scalar products in a Kulisch accumulator \textit{exact multiply add} (EMA). For an inner product, given a rounding function\footnote{$r(\cdot, b)$ is a rounding function that produces $b$ fractional bits, and $r_i(\cdot, b)$ is the $i$-th fractional bit returned. We assume IEEE 754-style round-to-nearest-even (with sticky bit OR-reduction) for $r(\cdot)$.} $r(\cdot)$ with the argument evaluated at infinite precision, EMA calculates $r(\sum_i a_ib_i)$, whereas FMA calculates $r(a_n b_n + r(a_{n-1} b_{n-1} + r(\cdots + r(a_1 b_1 + 0) \cdots)))$. Both EMA and FMA can be implemented for any floating point type. Gustafson proposed Kulisch accumulators to be standard for posits, terming them \textit{quires}.

Depending upon float dynamic range, EMA can be considerably more efficient than FMA in hardware. FMA must mutually align the addends $c$ and the product $ab$, including renormalization logic for subtraction cancellation, and the proper alignment cannot be computed until fairly late in the process. Extra machinery to reduce latency such as the \textit{leading zero (LZ) anticipator} or \textit{three path architectures} have been invented~\cite{Quinnell2007}. If multiply-add needs to be pipelined for timing closure, EMA knows upfront the location of the floating point of $c$ needed in alignment (as it is fixed), and can thus accumulate a new product into it every clock cycle, while a FMA must hold onto the starting value of the accumulator $c$ until later in the process, increasing the pipeline non-combinational area and often requiring greater use of an external register file (for multiple accumulators $c_i$ in concurrent use) and effective ``loop unrolling'' at software level to fill all pipeline slots. The rounding performed every FMA requires additional logic, and rounding error can still compound greatly across repeated sums.


\section{Multiplier efficiency}\label{sec:mul}

Floating point with EMA is still expensive, as there is added shifter, LZ counter, rounding, etc. logic. Integer MAC and float FMA/EMA both involve multiplication of fixed-point values; for int8/32 MAC this multiply is 63.4\% of the combinational power in our analysis at 28 nm (Section~\ref{sec:ASIC}).


A logarithmic number system (LNS)~\cite{kingsbury1971digital} avoids hardware multipliers entirely, where we round and encode $\log_B(x)$ for some base $B$ to represent a number $x \in \mathbb{R}$. Hitherto we have considered \textit{linear domain} representations, where $x \in \mathbb{R}$ is rounded and encoded as $x$ in integer, fixed or floating point representation (note that floating point is itself a combination of linear and log encodings).
Log domain operations on linear $x > 0, y > 0$ represented as $i = \log_2(x), j = \log_2(y)$ are:
\begin{align*}
\log_2(x \pm y) &= i + \sigma_{\pm}(j-i) \\
\log_2(xy) &= i + j  \mbox{\qquad\qquad\qquad \cite{coleman2000arithmetic}}\\
\log_2(x/y) &= i - j
\end{align*}
As values $x \leq 0$ are outside the log domain, sign and zero are handled separately~\cite{swartzlander1975sign}, as is $\pm \infty$. We encode $B=2$ log numbers with a sign bit and a signed fixed-point number of the form $m.f$, which represents the linear domain value $\pm 2^{(m + \sum_i f_i/2^i)}$.
For add/sub, without loss of generality, order $j \leq i$, and $\sigma_{\pm}(x) = \log_2(1 \pm 2^x)$; this is the historical weak point of a LNS, as implementations use costly LUTs or piecewise linear approximation of $\sigma_{\pm}(x)$. This can be more expensive than hardware multipliers. The approximation $\log_2(1 + x) \approx x$ for $x \in [0, 1]$ could also be used~\cite{miyashita2016convolutional}, but this adds significant error, especially with repeated sums.



$\sigma_{\pm}(x)$ need only be evaluated if one wishes to keep the partial sum in the log domain. As with Kulisch accumulation versus FMA, we accumulate in a different representation than the scalar product for efficiency. For $\sum_i a_i b_i$, we multiply $a_i b_i$ in the log domain, and then approximate as a linear domain floating point value for accumulation.
Translating log domain $m.f$ to linear is easier than $\sigma_{\pm}(x)$, as we can just consider the fractional portion $f$; $m$ is linear domain multiplication by $2^m$ (floating point exponent addition or fixed point bit shift). A LUT maps $f \in [0, 1)$ to $p(f) = 2^f - 1$. $p(f)$ is the linear representation of the log number fractional part; the LUT maps all bits of $f$ to a desired number of bits $\alpha$ of $p(f)$ or $r(p(f), \alpha)$, for a $(2^{f_{bits}} \times \alpha)$-bit LUT. Linear approximation of $m.f$ is the floating point value $\pm 2^m (1 + \sum_{i=1}^{\alpha} 2^{-i} r_i(p(f), \alpha))$. This is expanded in the usual way for Kulisch accumulation. Just as Kulisch accumulation is efficient for linear domain values up to a reasonably wide dynamic range, it proves quite efficient for our linear approximations of log values.

To convert a linear domain value back to log domain, we map $g \in [0, 1)$ to $q(g) = \log_2(1 + g)$. $g$ is a linear domain fixed-point fraction; to control the size of the LUT we only consider $\beta$ bits via rounding of $g$. $q(r(g, \beta))$ is similarly rounded to a desired $\gamma$ bits; note that this latter rounding is log domain. $r(q(r(g, \beta)), \gamma)$ is then a $(2^{\beta} \times \gamma)$-bit LUT. We also choose $\alpha \geq f_{bits} + 1, \beta \geq \alpha, \gamma = f_{bits}$ to ensure that log-to-linear-to-log conversion of $f$ is the identity, or $f = r(q(r(r(p(f), \alpha), \beta)), \gamma)$.

    We will name this (somewhat inaccurately) \textit{exact log-linear multiply-add} (ELMA). The log product and linear sum are each exact, but the log product is not represented exactly by $r(p(f))$ as this requires infinite precision, unlike EMA which is exact except for a final rounding. The intermediate log product avoids overflow or underflow with an extra bit for the product's $m$. If a linear-to-log mapping is desired (returning a log number after summation), there is also loss via $r(q(g))$.

    Combining log-to-linear mapping with Kulisch accumulation makes log domain multiply-add efficient and reasonably accurate. Small $p$ and $q$ LUTs reduce well in combinational logic. They are practical for 16-bit types too, as compression can be used to reduce the size. For larger types they are impractical, as $\alpha$, $\beta$, $\gamma$ need to scale with $2^{f_{bits}}$, at which point $\sigma_{\pm}$ is a better strategy. As with FMA, repeated summation via $\sigma_{\pm}$ is subject to magnitude difference error (\textit{e.g.}, the $c + ab = c$ case). Our approximation introduces error with $r(p(f))$ and $r(q(g))$, but mitigates repeated summation error and is immune to magnitude differences. This tradeoff seems acceptable in practice (Section~\ref{sec:analysis}).

An 8-bit log number by default suffers from the same problem as 8-bit IEEE-style floating point; the dynamic range is limited by the fixed point encoding. We can use the same tapering as used in $(N, s)$ posit for $m.f$ log numbers. $m$ is encoded as an exponent, and $f$ as a floating point significand. $f_{min}$ and $f_{max}$ are then exactly the same for posit-tapered base-2 log or linear domain values. Setting $\gamma = f_{bits}$ (which is at maximum $(N-3-s)$ for posits) introduces additional tapering rounding error, as subsequent rounding in encoding is performed outside regimes of maximum precision. $\gamma$ is increased up to 3 bits (guard, round and sticky bits in typical round-to-nearest-even) to improve accuracy here. This encoding we will refer to as $(N, s, \alpha, \beta, \gamma)$ log (posit tapered). We can similarly choose to encode log numbers using an IEEE 754 format (with biased exponents, NaN representations etc.); we use this for our ELMA comparison against float16 FMA in Section~\ref{sec:ASIC}.



\section {Additional hardware details}\label{sec:add-hw}

To make EMA/ELMA more energy efficient, we restrict accumulator range to $[f_{min}^2, f_{max}]$; handling temporary underflow rather than overflow is more important in our experience.
Kulisch accumulator conversion back to log or linear N-bit types uses a LZ counter and shifter but can be substantially amortized in two ways. First, many sums are performed, with final conversion done only once per inner product. Energy for the majority of work is thus lower than MAC/FMA (Section~\ref{sec:ASIC}); increased area for increased energy efficiency is generally useful in the era of ``dark silicon''~\cite{taylor2012dark}, or conversion module instances can be rationed (limiting throughput) and/or clock gated. Second, structures with local operand reuse (\textit{e.g.}, systolic arrays, fixed-function convolvers) naturally require fewer converter instances, reducing area (discussion in Section~\ref{sec:ASIC} as well). EMA and FMA accuracy are the same for a single sum $c + ab$; our power advantage would disappear in this domain, but the vast majority of flops/ops in NNs require repeated rather than singular sums. Note that int8/32 usage itself requires some conversion back to int8 in the end that we do not evaluate.

\section {FPGA experiments}\label{sec:analysis}

\begin{table}
{
\caption{\label{tab:resnet-acc}
ResNet-50 ImageNet validation set accuracy per math type
}}
\centering \begin{tabular}{llll}
\textbf{Math type} & \textbf{Multiply-add type} & \textbf{top-1 acc (\%)} & \textbf{top-5 acc (\%)} \\
\hline
float32 & FMA                  & 76.130                & 92.862 \\
\textbf{(8, 1, 5, 5, 7) log} & \textbf{ELMA}      & \textbf{-0.90} & \textbf{-0.20} \\[2mm]
(7, 1) posit & EMA            & -4.63 & -2.28 \\
(8, 0) posit & EMA            & -76.03 & -92.36 \\
(8, 1) posit & EMA            & -0.87 & -0.19 \\
(8, 2) posit & EMA            & -2.20 & -0.85 \\
(9, 1) posit & EMA            & -0.30 & -0.09 \\[2mm]
Jacob et al.~\cite{Jacob2018}: \\
float32 & FMA                  & 76.400 & n/a \\
int8/32 & MAC                  & -1.50  & n/a \\
Migacz~\cite{migacz}: \\
float32 & FMA                  & 73.230 \enspace       & 91.180 \\
int8/32 & MAC                 & -0.20  & -0.03\\


\end{tabular}
\medskip
\end{table}

Our implementation is in SystemVerilog for ASIC evaluation, built into an FPGA design with Intel FPGA OpenCL RTL integration support, with rudimentary PyTorch~\cite{paszke2017automatic} integration. Source code is available at \texttt{github.com/facebookresearch/deepfloat}.
We evaluate $(N, s)$ posit and $(N, s, \alpha, \beta, \gamma)$ log arithmetic on the ResNet-50 CNN~\cite{he2016deep} with the ImageNet ILSVRC12 validation set~\cite{russakovsky2015imagenet}. We use float32 trained parameters from the PyTorch model zoo, with batch normalization fused into preceding affine layers~\cite{Jacob2018}.
float32 parameters and network input are converted to our formats via round-to-nearest-even; no other adjustment of these values is performed.
When converting into or out of a Kulisch accumulator, we can add a small exponent bias factor, adjusting the input exponent by $m$, or the output exponent by $n$. This is effectively free (a small adder).
No changes are made to any activations except for such a bias of $n=-4$ at the last (fully connected) layer to recenter unnormalized log probabilities from around 16.0 to 1.0. Without this we have an additional loss in top-1 of around 0.5-1\%, with little change to top-5. If the Kulisch accumulator itself can be directly considered for top-$k$ comparison, this avoids the need as well. All math is replaced with the corresponding posit or log versions; average pooling is via division of the Kulisch accumulator.

Our results are in Table~\ref{tab:resnet-acc}, along with two int8/32 quantization comparisons. (8, 0) linear posit has insufficient dynamic range to work; activations are quickly rounded to zero. Our (8, 1, 5, 5, 7) log result remains very close to (8, 1) linear posit. The int8/32 results listed do not start from the same float32 parameters as our trained network, so they are not directly comparable. They use training with simulated quantization~\cite{Jacob2018} and KL-divergence calibration with sampled activations~\cite{migacz}, whereas we perform math in the usual way in our log or linear domain arithmetic after rounding input and parameters. We obtain reasonably similar precision without retraining, sampling activations or learning quantization parameters, while retaining general floating point representations in 8 bits.



\section {ASIC evaluation}\label{sec:ASIC}

We use Synopsys Design Compiler and PrimeTime PX with a commercially available 28 nm library, target clock 500 MHz. Process corners are SS@-40$^{\circ}$C synthesis, TT@25$^{\circ}$C power analysis at 0.81V. Table \ref{tab:pe} investigates multiply-add PEs, and as a proxy for an accelerator design, a 32x32 matrix multiplication systolic array with these PEs. The float16 FMA is Synopsys DesignWare dw\_fp\_mac. We accumulate to the C matrix in place (\textit{stationary C}), shifting out values upon completion. The int8/32 array outputs unprocessed int32; for ELMA, Kulisch accumulators are shifted across the PEs for C output and converted to 8 bit log at the boundary via 32 conversion/encoder modules. The 1024 PEs within do not include these (as discussed in Section~\ref{sec:add-hw}). 64 posit taper decoders are included for where A and B are passed as input. Power analysis uses testbench waves for 128-d vectors with elements drawn from $N(0, 1)$; int8 quantization has a max of $2\sigma$. PEs evaluate a variety of these inner products, and the systolic arrays a variety of GEMMs with these vectors.

\begin{table}
{
\caption{\label{tab:pe}
Chip area and power for 28 nm, 1-cycle multiply-add at 500 MHz
}}
\centering \begin{tabular}{lll}
\textbf{Component} & \textbf{Area} $\mu$m$^2$ & \textbf{Power} $\mu$W \\
\hline
\textbf{int8/32 MAC PE} & \textbf{336.672} & \textbf{283} \\
multiply & 121.212 & 108.0 \\
add & 117.810 & 62.3 \\
non-combinational & 96.768 & 112.7 \\[2mm]
\textbf{(8, 1, 5, 5, 7) log ELMA PE} & \textbf{376.110} & \textbf{272}\\
log multiply (9 bit adder) & 32.760 & 17.1 \\
$r(p(f))$ (16x5 bit LUT) & 8.946 & 5.4 \\
Kulisch shift (6 $\rightarrow$ 38 bit) & 81.774 & 71.0 \\
Kulisch add (38 bit) & 123.732 & 54.2 \\
non-combinational & 126.756 & 124.3 \\[2mm]
\textbf{float16 (w/o denormals) FMA PE} & \textbf{1545.012} & \textbf{1358} \\
\textbf{(5, 10) (11, 11, 10) log ELMA PE} & \textbf{1043.154} & \textbf{805} \\
(this log is (5, 10) float16-style encoding, same dynamic range;\\
denormals for log and float16 here are unhandled and flush to zero)\\[2mm]
\textbf{32x32 systolic w/ int8/32 MAC PEs} & \textbf{348231} & \textbf{226000} \\
\textbf{32x32 systolic w/ (8, 1, 5, 5, 7) log ELMA PEs } & \textbf{457738} & \textbf{195500} \\
\end{tabular}
\medskip
\end{table}

ELMA saves 90.9 $\mu$W over int8/32 on multiplication, but loses 68.3 $\mu$W on the add. ELMA non-combinational demands are higher with additional state required (Kulisch and decoded log numbers), but could be reduced by not handling underflow all the way to $f_{min}^2$. Despite the larger Kulisch adder, effectively only 6 bits are summed (with carry) each cycle versus up to 16 with int8/32; strategies for 500+ bit Kulisch accumulators~\cite{Uguen2017} might work in this small regime to further take advantage of this. Our 16-bit ELMA $\alpha=11$ $p(f)$ combinational LUT is 386 $\mu$m$^2$ despite compression, now a significant portion of the design. Larger $\alpha$ likely needs a compiled ROM or explicit compute of $p(f)$.

A more in-depth analysis for our work would need to determine a Pareto frontier between frequency/latency, per-operation energy, area, pipeline depth, math implementation and accuracy similar to the Galal et al. FPU generator work~\cite{galal2013fpu}, to see precisely in what regimes ELMA is advantageous. We provide our limited analysis here, however rough, to help motivate future investigation.

\section {Conclusions}

DNNs are resilient to many forms of numerical tinkering; they allow re-evaluation of design decisions made long ago at the bottom of the hardware stack with reduced fear of failure. The design space of hardware real number representations is indeed quite large and underexplored~\cite{Lindstrom2018}, as is as the opportunity to improve hardware efficiency and software simplicity with alternative designs and judicious use of numerical approximation. Log domain representations, posits, Kulisch accumulation and combinations such as ELMA show that floating point efficiency and applicability can be substantially improved upon. We plan on continuing investigation of this arithmetic design space at the hardware level with DNN training, and on general numerical algorithms in the future.

\subsubsection*{Acknowledgments}

We thank Synopsys for their permission to publish baseline and comparative results obtained by using their tools and DesignWare components, so we could present realistic numbers on our research using a popular 28 nm semiconductor technology node.

{\small
\bibliographystyle{ieee}
\bibliography{egbib}

\begin{thebibliography}{10}\itemsep=-1pt

\bibitem{cai2017deep}
Z.~Cai, X.~He, J.~Sun, and N.~Vasconcelos.
\newblock Deep learning with low precision by half-wave gaussian quantization.

\bibitem{coleman2000arithmetic}
J.~N. Coleman, E.~Chester, C.~I. Softley, and J.~Kadlec.
\newblock Arithmetic on the european logarithmic microprocessor.
\newblock {\em IEEE Transactions on Computers}, 49(7):702--715, 2000.

\bibitem{courbariaux2016binarized}
M.~Courbariaux, I.~Hubara, D.~Soudry, R.~El-Yaniv, and Y.~Bengio.
\newblock Binarized neural networks: Training deep neural networks with weights
  and activations constrained to+ 1 or-1.
\newblock {\em arXiv preprint arXiv:1602.02830}, 2016.

\bibitem{dettmers20158}
T.~Dettmers.
\newblock 8-bit approximations for parallelism in deep learning.
\newblock {\em arXiv preprint arXiv:1511.04561}, 2015.

\bibitem{DiCecco2017}
R.~DiCecco, L.~Sun, and P.~Chow.
\newblock Fpga-based training of convolutional neural networks with a reduced
  precision floating-point library.
\newblock In {\em 2017 International Conference on Field Programmable
  Technology (ICFPT)}, pages 239--242, Dec 2017.

\bibitem{fowers2018configurable}
J.~Fowers, K.~Ovtcharov, M.~Papamichael, T.~Massengill, M.~Liu, D.~Lo,
  S.~Alkalay, M.~Haselman, L.~Adams, M.~Ghandi, et~al.
\newblock A configurable cloud-scale dnn processor for real-time ai.
\newblock In {\em Proceedings of the 45th Annual International Symposium on
  Computer Architecture}, pages 1--14. IEEE Press, 2018.

\bibitem{galal2013fpu}
S.~Galal, O.~Shacham, J.~S. Brunhaver~II, J.~Pu, A.~Vassiliev, and M.~Horowitz.
\newblock Fpu generator for design space exploration.
\newblock In {\em Computer Arithmetic (ARITH), 2013 21st IEEE Symposium on},
  pages 25--34. IEEE, 2013.

\bibitem{golomb1966run}
S.~Golomb.
\newblock Run-length encodings (corresp.).
\newblock {\em IEEE transactions on information theory}, 12(3):399--401, 1966.

\bibitem{bfloat16}
Google.
\newblock {\em TPU TensorFlow ops}.
\newblock \url{https://cloud.google.com/tpu/docs/tensorflow-ops}.

\bibitem{gupta2015deep}
S.~Gupta, A.~Agrawal, K.~Gopalakrishnan, and P.~Narayanan.
\newblock Deep learning with limited numerical precision.
\newblock In {\em International Conference on Machine Learning}, pages
  1737--1746, 2015.

\bibitem{Gustafson2015}
J.~Gustafson.
\newblock {\em The End of Error: Unum Computing}.
\newblock Chapman \& Hall/CRC Computational Science. Taylor \& Francis, 2015.

\bibitem{Gustafson2017}
J.~L. Gustafson and I.~T. Yonemoto.
\newblock Beating floating point at its own game: Posit arithmetic.
\newblock {\em Supercomputing Frontiers and Innovations}, 4(2):71--86, 2017.

\bibitem{he2016deep}
K.~He, X.~Zhang, S.~Ren, and J.~Sun.
\newblock Deep residual learning for image recognition.
\newblock In {\em Proceedings of the IEEE conference on computer vision and
  pattern recognition}, pages 770--778, 2016.

\bibitem{horowitz20141}
M.~Horowitz.
\newblock 1.1 computing's energy problem (and what we can do about it).
\newblock In {\em Solid-State Circuits Conference Digest of Technical Papers
  (ISSCC), 2014 IEEE International}, pages 10--14. IEEE, 2014.

\bibitem{Jacob2018}
B.~Jacob, S.~Kligys, B.~Chen, M.~Zhu, M.~Tang, A.~Howard, H.~Adam, and
  D.~Kalenichenko.
\newblock Quantization and training of neural networks for efficient
  integer-arithmetic-only inference.
\newblock In {\em The IEEE Conference on Computer Vision and Pattern
  Recognition (CVPR)}, June 2018.

\bibitem{kingsbury1971digital}
N.~G. Kingsbury and P.~J. Rayner.
\newblock Digital filtering using logarithmic arithmetic.
\newblock {\em Electronics Letters}, 7(2):56--58, 1971.

\bibitem{koster2017flexpoint}
U.~K{\"o}ster, T.~Webb, X.~Wang, M.~Nassar, A.~K. Bansal, W.~Constable,
  O.~Elibol, S.~Gray, S.~Hall, L.~Hornof, et~al.
\newblock Flexpoint: An adaptive numerical format for efficient training of
  deep neural networks.
\newblock In {\em Advances in Neural Information Processing Systems}, pages
  1742--1752, 2017.

\bibitem{Kulisch2002}
U.~Kulisch.
\newblock {\em Advanced Arithmetic for the Digital Computer: Design of
  Arithmetic Units}.
\newblock Springer mathematics. Springer Vienna, 2002.

\bibitem{Kulisch2012}
U.~Kulisch.
\newblock {\em Computer Arithmetic and Validity: Theory, Implementation, and
  Applications}.
\newblock De Gruyter Studies in Mathematics. De Gruyter, 2012.

\bibitem{DBLP:journals/corr/LiL16}
F.~Li and B.~Liu.
\newblock Ternary weight networks.
\newblock {\em CoRR}, abs/1605.04711, 2016.

\bibitem{Lin2016}
D.~Lin, S.~Talathi, and S.~Annapureddy.
\newblock Fixed point quantization of deep convolutional networks.
\newblock In {\em International Conference on Machine Learning}, pages
  2849--2858, 2016.

\bibitem{Lindstrom2018}
P.~Lindstrom, S.~Lloyd, and J.~Hittinger.
\newblock Universal coding of the reals: alternatives to ieee floating point.
\newblock In {\em Proceedings of the Conference for Next Generation
  Arithmetic}, page~5. ACM, 2018.

\bibitem{migacz}
S.~Migacz.
\newblock 8-bit inference with tensorrt.
\newblock {\em Nvidia GTC}, 2017.

\bibitem{miyashita2016convolutional}
D.~Miyashita, E.~H. Lee, and B.~Murmann.
\newblock Convolutional neural networks using logarithmic data representation.
\newblock {\em arXiv preprint arXiv:1603.01025}, 2016.

\bibitem{morris1971tapered}
R.~Morris.
\newblock Tapered floating point: A new floating-point representation.
\newblock {\em IEEE Transactions on Computers}, 100(12):1578--1579, 1971.

\bibitem{muller2010handbook}
J.-M. Muller, F.~De~Dinechin, C.-P. Jeannerod, S.~Torres, et~al.
\newblock {\em Handbook of floating-point arithmetic}.
\newblock Springer, 2010.

\bibitem{paszke2017automatic}
A.~Paszke, S.~Gross, S.~Chintala, G.~Chanan, E.~Yang, Z.~DeVito, Z.~Lin,
  A.~Desmaison, L.~Antiga, and A.~Lerer.
\newblock Automatic differentiation in pytorch.
\newblock In {\em NIPS-W}, 2017.

\bibitem{Quinnell2007}
E.~Quinnell, E.~E. Swartzlander, and C.~Lemonds.
\newblock Floating-point fused multiply-add architectures.
\newblock In {\em Signals, Systems and Computers, 2007. ACSSC 2007. Conference
  Record of the Forty-First Asilomar Conference on}, pages 331--337. IEEE,
  2007.

\bibitem{rastegari2016xnor}
M.~Rastegari, V.~Ordonez, J.~Redmon, and A.~Farhadi.
\newblock Xnor-net: Imagenet classification using binary convolutional neural
  networks.
\newblock In {\em European Conference on Computer Vision}, pages 525--542.
  Springer, 2016.

\bibitem{russakovsky2015imagenet}
O.~Russakovsky, J.~Deng, H.~Su, J.~Krause, S.~Satheesh, S.~Ma, Z.~Huang,
  A.~Karpathy, A.~Khosla, M.~Bernstein, et~al.
\newblock Imagenet large scale visual recognition challenge.
\newblock {\em International Journal of Computer Vision}, 115(3):211--252,
  2015.

\bibitem{swartzlander1975sign}
E.~E. Swartzlander and A.~G. Alexopoulos.
\newblock The sign/logarithm number system.
\newblock {\em IEEE Transactions on Computers}, 100(12):1238--1242, 1975.

\bibitem{taylor2012dark}
M.~B. Taylor.
\newblock Is dark silicon useful? harnessing the four horsemen of the coming
  dark silicon apocalypse.
\newblock In {\em Design Automation Conference (DAC), 2012 49th ACM/EDAC/IEEE},
  pages 1131--1136. IEEE, 2012.

\bibitem{DBLP:journals/corr/abs-1712-02427}
A.~Tulloch and Y.~Jia.
\newblock High performance ultra-low-precision convolutions on mobile devices.
\newblock {\em CoRR}, abs/1712.02427, 2017.

\bibitem{Uguen2017}
Y.~Uguen and F.~De~Dinechin.
\newblock Design-space exploration for the kulisch accumulator.
\newblock 2017.

\bibitem{vanhoucke2011improving}
V.~Vanhoucke, A.~Senior, and M.~Z. Mao.
\newblock Improving the speed of neural networks on cpus.
\newblock Citeseer.

\bibitem{wilkinson1963rounding}
J.~H. Wilkinson.
\newblock {\em Rounding errors in algebraic processes}.
\newblock Prentice-Hall, 1963.

\bibitem{zuras2008ieee}
D.~Zuras, M.~Cowlishaw, A.~Aiken, M.~Applegate, D.~Bailey, S.~Bass,
  D.~Bhandarkar, M.~Bhat, D.~Bindel, S.~Boldo, et~al.
\newblock Ieee standard for floating-point arithmetic.
\newblock {\em IEEE Std 754-2008}, pages 1--70, 2008.

\end{thebibliography}
}

\end{document}